# A Distributed Gradient Approach for System Optimal Dynamic Traffic Assignment

Mehrzad Mehrabipour and Ali Hajbabaie[1]

*Abstract*—This study presents a distributed gradient-based approach to solve system optimal dynamic traffic assignment (SODTA) formulated based on the cell transmission model. The algorithm distributes SODTA into local sub-problems, who find optimal values for their decision variables within an intersection. Each sub-problem communicates with its immediate neighbors to reach a consensus on the values of common decision variables. A sub-problem receives proposed values for common decision variables from all adjacent sub-problems and incorporates them into its own offered values by weighted averaging and enforcing a gradient step to minimize its objective function. Then, the updated values are projected onto the feasible region of the sub-problems. The algorithm finds high quality solutions in all tested scenarios with a finite number of iterations. The algorithm is tested on a case study network under different demand levels and finds solutions with at most a 5% optimality gap.

## I. INTRODUCTION

Dynamic Traffic Assignment (DTA) is a well-studied research area to determine time-dependent traffic flows by minimizing the cost of the system or individual users in transportation networks. DTA deployment has brought many benefits to a wide range of applications over the past decades. Network design [1]–[4], traffic operations [5], [6], congestion pricing [7], evacuation planning [8], and traffic management systems [9] are some of DTA applications.

DTA with accurate network loading models has many decision variables and constraints to encompass its spatial/temporal scales and the number of origin-destination (OD) pairs. The number of decision variables and constraints can further increase based on the utilized network loading concept. For instance, a cell transmission model (CTM)-based DTA will have more decision variables and constraints than a link performance function-based DTA, as CTM divides each link into several shorter links, which increases the number of decision variables.

Centralized approaches solve the optimization model with no decomposition, distribution, or parallelism. Central optimization frameworks do not scale with the size of DTA, especially when complex network loading concepts, such as CTM are in use [5], [10]–[13]. Decomposition approaches convert the model into sub-problems and a master problem. Each subproblem has less computational complexity in comparison with the original model because the subproblem contains only a fraction of the variables and constraints of the original model. Therefore, decomposition approaches can improve scalability; however, non-spatial decomposition approaches face computational complexity growth in the sub-problems as the size of the original problem grows. Besides, all existing decomposition techniques require solving the master problem (a central component) to aggregate the solutions found by sub-problems and expanding the network will increase the complexity of the master problem as well. Therefore, while decomposition approaches can provide optimality bounds and scale better than the centralized approaches, they eventually become intractable when the size of the problem grows [14]–[19].

This paper presents a distributed gradient-based approach (DGA) to overcome the discussed drawbacks of existing decomposition approaches. The proposed methodology distributes network-level DTA into several intersection-level DTA sub-problems. Constraints and the objective function are exclusively distributed among sub-problems.

It allocates a computational node to each sub-problem and does not require a central component or a master problem. Therefore, it scales better than the decomposition algorithms, whose complexity depend on the network size. The computational complexity of sub-problems does not depend on the number of nodes and links in this approach as a result of the spatial decomposition; however, it is a function of the number of OD pairs. This intersection-level distribution is well-suited for urban network planning purposes because extending network by adding more intersections (nodes and links) will not change the computational complexity of the methodology and its structure.

The approach starts with setting initial values for decision variables. These values can be found by simulating the network with a path-based CTM and sending the demand through the shortest paths found by Dijkstra's algorithm. Then, the approach performs three steps to update the value of decision variables iteratively at each sub-problem: (1) it first incorporates the value of common decision variables among sub-problems by taking a weighted average, (2) the approach adjusts the value of the decision variables from the first step by moving them towards the negative direction of the gradient of the objective function at each sub-problem to minimize the objective, and (3) it projects the values of step 2 onto the set of constraints of each sub-problem to maintain feasibility. The approach iterates among these three steps until the values of shared decision variables from sub-problems reach consensus: there is an acceptable disagreement on the generated values of the decision variables. The sub-problems communicate based on a hypothetical information exchange graph that consists of nodes and directional links. Each node and link represent a sub-problem and the communication between sub-problems, respectively. A weight is assigned to each link for incorporating information among sub-problems. The approach uses these weights for taking a weighted average in the first step.

[1] Corresponding Author: Ali Hajbabaie, e-mail: ahajbab@ncsu.edu. The authors are with the Department of Civil, Construction, and Environmental Engineering at North Carolina State University.

It can be shown theoretically that the algorithm eventually converges to the optimal solution of the problem. Note that the approach can work without necessarily starting with a feasible solution. Furthermore, the approach is applicable to problems with nonlinear and quadratic objective functions since it needs either the gradient or sub-gradients of the objective function.

In the remainder of this paper, we present a review of relevant literature to solve SODTA problems. Then, we describe the problem formulation and the solution technique in separate sections. The discussion is continued by presenting a case study network and numerical results. Finally, the concluding remarks and directions for future research are presented.

## II. LITERATURE REVIEW

Ziliaskopoulos [10] formulated SODTA with one destination using CTM as a linear and convex program. He solved the program centrally using the Simplex method for a small network with 10 cells. Beard and Ziliaskopoulos [5] further improved this formulation to jointly optimize signal timing parameters and traffic assignment considering multiple OD pairs. This formulation was a mixed-integer linear program and was solved centrally for only two intersections. Chiu and Zheng [11] also modified the CTM-based SODTA to address emergency and disaster conditions. They presented a linear formulation with multiple OD pairs and optimized it for a simplified network of 40 cells using the interior point method to find paths and departure schedules for muti-priority groups. Aziz and Ukkusuri [13] proposed a non-linear and quadratic CTM-based SODTA formulation to minimize emission and travel time with a single destination. They solved the model centrally for a network with 88 cells. Zheng and Chiu [12] proposed a network flow algorithm to solve Ziliaskopoulos's [10] formulation and overcome the scalability issue of the previous centralized approaches. They first proved that SODTA was equivalent to the earliest arrival flow problem. Then, they created a time-expanded network to solve the earliest arrival flow problem. They applied the algorithm to a network with 2012 cells. Even though the algorithm was more efficient than the previous techniques, it was restricted to single destination problems since the earliest arrival flow problem was not equivalent to SODTA for multiple OD pairs. Therefore, the algorithm was not applicable to SODTA with multiple OD pairs.  All these approaches used a centralized solution approach, which did not scale well with temporal and spatial scales of the CTM-based traffic assignment problems. Therefore, they were appicable to smaller-size problems with a limited amount of available memory and time.

Long et al. [20] showed that the link-based formulation was more computationally efficient compared to path-based counterparts. However, solving the problem for larger cases still required more efficient algorithms. They presented various formulations to model SODTA using Link Transmission Mode (LTM) to optimize route and departure choice decisions. Different combinations of models including first-in-first-out (FIFO) and non-holding-back constraints were presented. The problems without FIFO constraints were solved using a commercial solver and those with FIFO constraints were optimized using the branch and bound algorithm for a network of twelve nodes. Long and Szeto [21] also studied LTM for the SODTA problem considering FIFO conditions.  Vehicle holding back was also addressed using mixed-integer constraints. The proposed formulation with FIFO constraints was a mixed-integer linear program that had more computational complexity compared to linear problems. The problem was solved using a branch and bound algorithm via commercial software.

Decomposition approaches were developed to decrease the computational complexity of transportation problems [22]–[27] and SODTA. Li et al. [17] presented a decomposition algorithm based on the Dantzig-Wolfe principle to solve Ziliaskopoulos's [10] formulation. The algorithm consisted of a sub-problem and a master problem. The sub-problem generated link flows using a minimum-cost-flow problem, and the master problem found the best convex combination of all generated flows. The number of decision variables in the master problem increases over iterations, and the computational complexity of the sub-problems and the master problem increase by adding more nodes and links to the network. Moreover, the algorithm cannot be applied to cases where the ratio of free-flow speed to backward propagation speed is not an integer value. Lin et al. [16] also used the Dantzig-Wolfe principle to develop a heuristic-based decomposition algorithm to calibrate the flow capacity of a User Equilibrium (UE) DTA with a single destination and CTM traffic dynamics. They created sub-problems by decomposing the dual formulation for bi-level capacity calibration problem. The master problem considered the relaxed constraints in sub-problems and combined their solutions. Similar to Li et al.'s [17] study, the growth in the computational complexity of the master problem and sub-problems is unavoidable. Ramadurai and Ukkusuri [28] proposed an approach with acyclic sub-networks for ODs and activity sequences to solve a user equilibrium DTA problem with a point queue model. The route choice, location, time, and duration of an activity were decision variables. In the approach, sub-network equilibrium was achieved by changing the flow of paths. The optimality condition was not guaranteed in this simulation-based approach.

Jafari et al. [15] decomposed the UE static traffic assignment problem spatially by creating sub-networks to represent sub-problems. The sub-problems found the traffic assignment solution within a sub-network given the regional demand from a master problem. The master problem solved the static traffic assignment problem for an aggregated network using modified travel times obtained from the sub-problems. They needed to create an aggregated network with arbitrary arcs that contained all paths based on the original network. They applied a method of sensitivity analysis by Jafari and Boyles [29] to find required parameters for arbitrary links and relative changes in path flows by demand fluctuations at each iteration. The approach by Boyles [30] provided the basis for network aggregation to create the master problem and model sub-networks. Boyles's [30] method could find the cost function of arbitrary arcs in the aggregated network as well. In Jafari et al.'s [15] approach, considering more sub-networks leads to increasing the complexity of the aggregated network and the master problem. Moreover, the algorithm is restricted to acyclic networks and requires an initial feasible solution that satisfies UE conditions.

Mehrabipour et al. [14] developed an OD-based decomposition approach using the Dantzig-Wolfe principle to overcome some drawbacks in previous approaches. The approach solved a CTM-based SODTA problem and consisted of a master problem and sub-problems. Each sub-problem found cell flows for one OD pair, and the master problem combined the generated flows considering the unseen constraints in sub-problems on total flow over all OD pairs. This approach increased the scalability of centralized approaches while guaranteeing optimality conditions. However, the computational complexity of the master problem grew with the number of iterations. Besides, the computational complexity of the master and sub-problems increased by adding more nodes and links to the network.

*A. The Summary of Contributions*

The contributions of this paper are as follow:
1. The proposed approach is completely distributed. The sub-problems work cooperatively without requiring a central optimization unlike existing decomposition approaches developed in [14]–[19], [31]–[33] Therefore, the proposed approach is scalable with the size of the network.
2. The approach finds high quality solutions to the SODTA problem with a reasonable number of iterations . The required assumptions for convergence do not restrict the application of the approach to specific network properties such as limiting the ratio of free-flow speed to backward propagation wave speed [17] and having specific network geometry [15], [19].
3. The computational complexity of sub-problems is independent of the number of nodes and links in the network since we create each sub-problem by a spatial distribution of the objective function and the set of constraints unlike the previous approaches: (e.g., [15]–[19], [31]–[33]).

III. PROBLEM FORMULATION

SODTA is formulated as a linear program utilizing the CTM traffic dynamics introduced by Daganzo [34], [35]. In cell transmission model, a road is discretized into homogenous segments called cells. Homogenous cells have the same length and free-flow speed. Also, the study period is divided into time intervals named time steps. The length of each cell is equivalent to the distance that a vehicle can travel in one-time step. Note that CTM relates flow and density in each cell using non-linear equations. Some examples of formulations based on CTM can be found in [36]–[39]. We used the linearized form of the CTM-based SODTA formulation introduced by Beard and Ziliaskopoulos [5] and modified the set of OD pairs to reduce the computational complexity. The cell-based constraints developed in [40] can be incorporated in our formulation to resolve priority issues at merge cells in the future.

TABLE VI (see Appendix) presents the sets, decision variables, parameters, notation, and terms used throughout the paper. Let $C, T, C_{OD}$, and $S(i)$ respectively denote the set of cells, time steps, OD pairs, and cells successor to cell $i \in C$.

This formulation includes two sets of decision variables. The first set is the number of vehicles $x_i^{t,od}$ in cell $i \in C$ at time step $t \in T$ associated with OD pair $(o,d) \in C_{OD}$, and the second set is the number of vehicles $y_{ij}^{t,od}$ flowing from cell $i \in C$ to successor cell $j \in S(i)$ at time step $t \in T$ associated with OD pair $(o,d) \in C_{OD}$. Traffic signal timing parameters are input to this formulation, and equation (1) finds the variable saturation flow rate $f_i^t$ at intersection cell $i \in C_I$ for time step $t \in T$ based on the signal timing parameters. We use $C_I$ and $F_i$ to denote the set of intersection cells and constant saturation flow rate of cell $i \in C_I$. The signal status $g_i^t$ is a binary parameter defined for all intersection cells $i \in C_I$ and time steps $t \in T$. When the signal is green, it will be one, and zero otherwise. The variable saturation flow rate is equal to the constant saturation flow rate if the signal in the intersection cell is green.

$$f_i^t = g_i^t F_i \qquad \forall\, t \in T, i \in C_I \qquad (1)$$

We define the minimization of total travel time as the objective function in expression (2). The total travel time is found by summing the number of all vehicles $x_i^{t,od}$ in all network cells except for sink cells $i \in C \setminus C_S$ for all OD pairs $(o,d) \in C_{OD}$ over all time steps $t \in T$ and multiplying the result by $\tau$ (i.e., the duration of each time step). We can eliminate the time step duration $\tau$ from the objective function since it has a constant value.

$$\text{Min } \sum_{(o,d)\in C_{OD}} \sum_{t \in T} \sum_{i \in C \setminus C_S} \tau x_i^{t,od} \qquad (2)$$

Constraints (3), (4), and (5) show the conservation of flow for ordinary, source, and sink cells, respectively. An increase or decrease in the number of vehicles $x_i^{t+1,od} - x_i^{t,od}$ between time steps $t \in T$ and $t+1 \in T$ is equal to the difference of the total inflow and outflow of cell $i \in C$ at time step $t \in T$ for OD pair $(o,d) \in C_{OD}$.

$$\sum_{k \in P(i)} y_{ki}^{t,od} - \sum_{j \in S(i)} y_{ij}^{t,od} = x_i^{t+1,od} - x_i^{t,od} \qquad (3)$$
$$\forall\, t \in T, i \in C \setminus \{C_S, C_O\}, (o,d) \in C_{OD}$$

$$\mathcal{D}_i^{t,od} - \sum_{j \in S(i)} y_{ij}^{t,od} = x_i^{t+1,od} - x_i^{t,od} \qquad (4)$$
$$\forall\, t \in T, i \in C_O, (o,d) \in C_{OD}$$

$$\sum_{i \in P(j)} y_{ij}^{t,od} = x_j^{t+1,od} - x_j^{t,od} \qquad (5)$$
$$\forall\, t \in T, j \in C_S, (o,d) \in C_{OD}$$

Constraint (6) limits the total outgoing flow $\sum_{j \in S(i)} y_{ij}^{t,od}$ from cell $i \in C$ to its successor cells $j \in S(i)$ to the occupancy $x_i^{t,od}$ of the cell at time step $t \in T$ with OD pair $(o,d) \in C_{OD}$.

$$\sum_{j \in S(i)} y_{ij}^{t,od} \leq x_i^{t,od} \qquad \forall\, t \in T, i \in C, (o,d) \in C_{OD} \qquad (6)$$

Constraints (7) and (8) respectively ensure that the total outgoing flow from and the total incoming flow to a cell are limited to the constant saturation flow rate of the cell. We use $F_i$ to denote the constant saturation flow rate of cell $i \in C$.

$$\sum_{(o,d)\in C_{OD}} \sum_{j \in S(i)} y_{ij}^{t,od} \leq F_i \qquad \forall\, t \in T, i \in C \qquad (7)$$

$$\sum_{(o,d)\in C_{OD}} \sum_{i \in P(j)} y_{ij}^{t,od} \leq F_j \qquad \forall\, t \in T, j \in C \qquad (8)$$

The total incoming flow $\sum_{(o,d)\in C_{OD}} \sum_{i \in P(j)} y_{ij}^{t,od}$ to cell $j \in C$ at time step $t \in T$ should be less than or equal to the available capacity $M_j - \sum_{(o,d)\in C_{OD}} x_j^{t,od}$ of that cell as shown by constraint (9). Note that $M_j$ denotes the maximum number of vehicles that cell $j \in C$ can accommodate. We use $\delta$ to denote the ratio of free-flow speed to the backward propagation speed.

$$\sum_{(o,d)\in C_{OD}} \sum_{i\in P(j)} y_{ij}^{t,od} \leq \delta\left(M_j - \sum_{(o,d)\in C_{OD}} x_j^{t,od}\right) \quad (9)$$

$$\forall\, t \in T, j \in C$$

Constraint (10) limits the total outgoing flow $\sum_{(o,d)\in C_{OD}} \sum_{j\in S(i)} y_{ij}^{t,od}$ from intersection cell $i \in C_I$ to variable saturation flow rate $f_i^t$ at time step $t \in T$. We sum the flow over all links between cell $i \in C_I$ and its successor cell $j \in S(i)$ over all OD pairs to find the total outflow $\sum_{(o,d)\in C_{OD}} \sum_{j\in S(i)} y_{ij}^{t,od}$. Equation (1) finds the variable saturation flow rate $f_i^t$ given the signal status $g_i^t$ for intersection cell $i \in C_I$ and time step $t \in T$.

$$\sum_{(o,d)\in C_{OD}} \sum_{j\in S(i)} y_{ij}^{t,od} \leq f_i^t \quad \forall\, t \in T, i \in C_I \quad (10)$$

Constraints (11) and (12) are used to ensure that the decision variables are nonnegative.

$$x_i^{t,od} \geq 0 \quad \forall\, t \in T, i \in C, (o,d) \in C_{OD} \quad (11)$$
$$y_{ij}^{t,od} \geq 0 \quad \forall\, t \in T, i \in C\,\{C_S\}, j \in S(i), (o,d) \in C_{OD} \quad (12)$$

The summary of the proposed optimization model follows:
Min $\sum_{(o,d)\in C_{OD}} \sum_{t\in T} \sum_{i\in C \setminus C_S} \tau x_i^{t,od}$
s.t.
(3)-(12)

## IV. METHODOLOGY

We present a distributed gradient-based methodology to solve the SODTA problem. The discussions are continued in five subsections: (1) the distribution of SODTA problem formulation, (2) initialization, (3) the procedure of distributed gradient-based update, (4) termination criteria, and (5) convergence properties.

In the formulation distribution section, we partition the cell-based SODTA formulation among sub-problems with an intersection-level segmentation so that the number of sub-problems is equal to the number of intersections. Each sub-problem contains some parts of the objective function and constraints of the original SODTA formulation that have decision variables corresponding to the cells and links within the intersection assigned to the sub-problem. The sub-problems do not share any constraint, and the objective function is distributed among sub-problems so that the summation of the objective functions is equivalent to the objective function of the original SODTA problem.

In the initialization step, we find initial values for the decision variables of each sub-problem. We start by generating the shortest paths for each OD pair using Dijkstra's algorithm [41] and sending the demand to the network through these paths using a path-based CTM simulation [14], [42]. We initialize the approach with the occupancy and flow values that are the outputs of the simulation. This simulation is a cell transmission model developed by Daganzo (1990) and used to simulate the movement of vehicles across the network. Then, the value of the decision variables of each sub-problem is updated over iterations. The updated value for a variable will be the weighted average of all presented values for that variable. We find the weighted average at the current iteration for all sub-problems using the value of decision variables from either the initialization step or the previous iteration. Then, the values of variables are updated considering the negative direction of the gradient of the objective function of each sub-problem (to minimize it) and project the value of the decision variables on the set of constraints at each sub-problem to maintain feasibility. The approach iterates until the conflict among the proposed values from the sub-problems is within an acceptable threshold. Note that having a unique optimal solution is not guaranteed in our approach. Fig. 1 shows the overall framework of the methodology, and we provide more details about each part of the figure in the rest of this section.

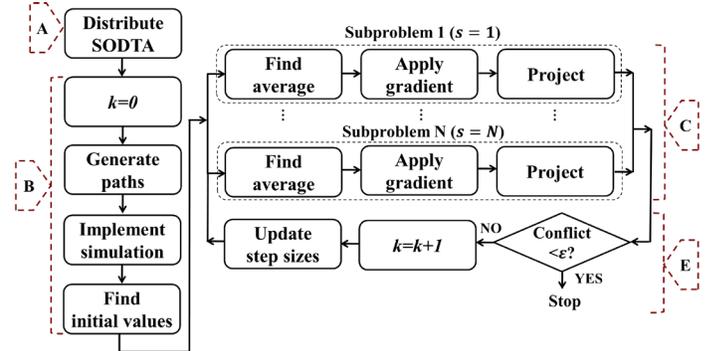

Fig. 1. The framework of the distributed gradient-based methodology

### B. The Distribution of SODTA Formulation

The first step of the methodology is to decompose the network-level SODTA problem formulation into intersection-level sub-problems. Each sub-problem will have decision variables associated with one intersection. The summation of the objective functions from each sub-problem is equivalent to the objective function of the original SODTA problem, and the union of the constraints from each sub-problem is equivalent to the original constraint set. Note that the sub-problems do not share any constraints. This intersection-level decomposition is well-suited for solving problems in urban street networks since assigning more intersections to the network will not change the architecture of the methodology and its computational complexity. Regional partitioning is used in [7] for static traffic assignment and in [43] for SODTA with route and departure time variables. More general regional partitioning examples can be found in [44]–[46]. The boundary of each intersection should be selected such that each intersection (subproblem) has more or less the same number of decision variables to balance out the run time and reduce overhead delays. A general formulation for sub-problem $s \in N$ is shown below. This formulation differs from the original formulation in the set of cells for which the constraints are defined. The sets are defined in Table VI. A formulation for sub-problem $s \in N$ is shown in (13)-(23).

$$\text{Min } \sum_{(o,d)\in C_{OD}} \sum_{t\in T} \sum_{i\in Ces^s} x_i^{t,od} \quad (13)$$
$$\sum_{k\in P(i)} y_{ki}^{t,od} - \sum_{j\in S(i)} y_{ij}^{t,od} = x_i^{t+1,od} - x_i^{t,od} \quad (14)$$
$$\forall\, t \in T, i \in Ceos^s, (o,d) \in C_{OD}$$
$$\mathcal{D}_i^{t,od} - \sum_{j\in S(i)} y_{ij}^{t,od} = x_i^{t+1,od} - x_i^{t,od} \quad (15)$$
$$\forall\, t \in T, i \in Co^s, (o,d) \in C_{OD}$$
$$\sum_{i\in P(j)} y_{ij}^{t,od} = x_j^{t+1,od} - x_j^{t,od} \quad (16)$$

$$\sum_{j\in S(i)} y_{ij}^{t,od} \leq x_i^{t,od} \quad \forall\, t \in T, j \in Cs^s, (o,d) \in C_{OD} \quad (17)$$

$$\sum_{(o,d)\in C_{OD}} \sum_{j\in S(i)} y_{ij}^{t,od} \leq F_i \quad \forall\, t \in T, i \in Ces^s, (o,d) \in C_{OD} \quad (18)$$

$$\sum_{(o,d)\in C_{OD}} \sum_{i\in P(j)} y_{ij}^{t,od} \leq F_j \quad \forall\, t \in T, j \in Ceo^s \quad (19)$$

$$\sum_{(o,d)\in C_{OD}} \sum_{i\in P(j)} y_{ij}^{t,od} \leq \delta(M_j - \sum_{\forall (o,d)\in C_{OD}} x_j^{t,od}) \quad (20)$$

$$\forall\, t \in T, j \in Ceo^s$$

$$\sum_{(o,d)\in C_{OD}} \sum_{j\in S(i)} y_{ij}^{t,od} \leq f_i^t \quad \forall\, t \in T, i \in Ci^s \quad (21)$$

$$x_i^{t,od} \geq 0 \quad \forall\, t \in T, i \in C^s, (o,d) \in C_{OD} \quad (22)$$

$$y_{ij}^{t,od} \geq 0 \quad \forall\, t \in T, i \in C^s, j \in S(i), (o,d) \in C_{OD} \quad (23)$$

*C. Initialization*

The approach starts with initial values for all decision variables at the first iteration, see Fig. 1. The initial solutions do not have to be feasible for the original SODTA formulation. We first implement Dijkstra's algorithm [41] to generate the shortest paths for all OD pairs. Note that we assume that the network is under free-flow condition. Then, we use the path-based CTM simulation introduced by Ukkusuri et al. [42] to find occupancy $x_i^{t,od}$ for cell $i \in C^s$ at time step $t \in T$ with OD pair $(o,d) \in C_{OD}$ and flow $y_{ij}^{t,od}$ for links between cell $i \in C^s$ and its successor cell $j \in S(i)$ at time step $t \in T$ with OD pair $(o,d) \in C_{OD}$ for all sub-problems $s \in N$.

*D. Distributed Gradient Update*

This procedure updates the decision variables of all sub-problems in iteration $k+1 \in K$ using three main steps:

*Step 1*: Each sub-problem optimizes the values of its decision variables. Therefore, decision variables that are in common between several sub-problems will have various values. This step sets the value of these decision variables in each sub-problem equal to their weighted average. The weights in each sub-problem are determined such that they satisfy the required conditions for convergence.

*Step 2*: The approach moves the computed values for sub-problem $s \in N$ from *Step 1* towards the negative direction of the gradient of the objective function of the corresponding sub-problem to minimize the total travel time within the region assigned to the sub-problem.

*Step 3*: The approach projects the decision variable values in sub-problem $s \in N$ from *Step 2* onto the set of constraints of the sub-problem to make the values feasible for that sub-problem.

We describe these steps with a mathematical representation in the rest of this section. We first introduce four definitions for the information exchange graph (Definition 1), neighbors of a sub-problem (Definition 2), the gradient of a function (Definition 3), and the projection operator (Definition 4).

**Definition 1** (Information exchange graph): The information exchange graph $G = (N, A)$ contains nodes and directional links that belong to sets $N$ and $A$, respectively. Node $s \in N$ represents a sub-problem, and link $(s, \acute{s}) \in A$ shows the transfer of information from sub-problem $s \in N$ to sub-problem $\acute{s} \in N$. If there is a variable in common between sub-problems $s, \acute{s} \in N$, directional links $(s, \acute{s}) \in A$ and $(\acute{s}, s) \in A$ are needed. There is also a self-arc at each node i.e., $\{(s, \acute{s}): s = \acute{s}, \forall s, \acute{s} \in N\}$.

This arc represents the use of information in the computation of sub-problem $s \in N$ that is generated by itself. Weight value $w_{s\acute{s}}$ for the link going from node $s \in N$ to node $\acute{s} \in N$ is also assigned to link $(s, \acute{s}) \in A$ to be used for incorporating information from sub-problem $s \in N$ for decision variables in sub-problem $\acute{s} \in N$. Assumption 1 determines the value of weights and is required to prove the convergence of the approach. We can use identical values for wights on all links $(s, \acute{s}) \in A$ if $s \neq \acute{s}$. Weights on self-arcs $s \in N$ satisfy equation $w_{ss} = -\sum_{\acute{s}\in C^{\mathfrak{n}^s}:s\neq \acute{s}} w_{\acute{s}s}$. Numerical examples for the weights are provided for the simple example in this section and the test network in the result section.

**Assumption 1** Let us use $\theta$ and $\theta'$ to show a lower and an upper bound for weights $w_{s\acute{s}}$ on all links $(s, \acute{s}) \in N$, respectively. The information exchange graph $G = (N, A)$ and weight $w_{s\acute{s}}$ on link $(s, \acute{s}) \in N$ should satisfy the following conditions:

(1) Information exchange graph $G = (N, A)$ should be connected.
(2) Weights must satisfy $w_{s\acute{s}} = w_{\acute{s}s}$ for all $s, \acute{s} \in N$
(3) If weight $w_{s\acute{s}} \neq 0$, we have $\theta \leq w_{s\acute{s}} \leq \theta'$ for all $s, \acute{s} \in N$ such that $\theta, \theta' > 0$
(4) Weights on self-arcs $s \in N$ must satisfy equation $w_{ss} = -\sum_{\acute{s}\in C^{\mathfrak{n}^s}:s\neq \acute{s}} w_{\acute{s}s}$

**Definition 2** (The neighbors of a sub-problem): The neighbors of sub-problem $s \in N$ are sub-problems $\acute{s} \in N: \acute{s} \neq s$ that offer estimation for at least one decision variable of sub-problem $s \in N$. In other words, sub-problem $\acute{s} \in N$ is a neighbor of sub-problem $s \in N$ if and only if $w_{\acute{s}s} > 0$. We define the set of neighbors of sub-problem $s \in N$ including itself by $\mathfrak{n}^s$.

In general, any node $s \in N$ may be connected to any other node $\acute{s} \in N$ in the information exchange graph $G$ even if the sub-problems (nodes) are not immediate neighbors in the original (physical) network if their corresponding sub-problems share a decision variable. This fact will not affect the distributed structure of the methodology because the approach uses the information from the previous iteration for the exchange process, not the current iteration. However, the structure of SODTA formulation and intersection-based distribution lead to the presence of links only between immediate neighbors in the information exchange graph. The reason is that the immediate neighbors share decision variables corresponding to the links between any two intersections (regions).

**Definition 3** (Gradient of a function): Let $\mathcal{F}_s(\boldsymbol{x}), \mathcal{G}_s \in \mathbb{R}^n$, and $X_s$ respectively denote the objective function value of sub-problem $s \in N$ given vector $\boldsymbol{x}$, the gradient of the objective function $\mathcal{F}_s(\boldsymbol{x})$, and the feasible region of sub-problem $s \in N$. The gradient $\mathcal{G}_s$ satisfies inequality (24) for all vectors $\boldsymbol{z}, \boldsymbol{x} \in X_s$.

$$\mathcal{F}_s(\boldsymbol{z}) + \mathcal{G}_s'(\boldsymbol{x} - \boldsymbol{z}) \leq \mathcal{F}_s(\boldsymbol{x}) \quad (24)$$

**Definition 4** (Projection operator): We use the projection operator $\boldsymbol{P}_X[\boldsymbol{z}]$ to find the projection of vector $\boldsymbol{z}$ onto a closed convex set $X$ using Euclidean norm as shown in (25).

$$\boldsymbol{P}_X[\boldsymbol{z}] = argmin_{x\in X} \|\boldsymbol{z} - \boldsymbol{x}\| \quad (25)$$

A projection is a linear transformation, and the projection operator is used to map any vector onto a closed convex set. By solving $\arg\min_{z \in X}\|z - x\|$, we can apply this operator to map vector of $x \in X$ on set $X$.

We now describe all three steps for iteration $k + 1 \in K$ assuming that the value of occupancy and flow decision variables are available (either from the initialization step or previous iteration $k \in K$ for all sub-problems $s \in N$). We first update the value of decision variables as described in *Steps 1 and 2*. We define auxiliary parameters $p_i^{t,od}: t \in T, i \in C^s, (o,d) \in C_{OD}$ and $q_{ij}^{t,od}: t \in T, i \in C^s, j \in S(i), (o,d) \in C_{OD}$ for updating the value of cell occupancy and flow decision variables, respectively. If a decision variable appears in only one sub-problem, we find its auxiliary parameter using either equation (26) or (27). Equations (26) and (27) find the auxiliary parameters for occupancy and flow variables, respectively. Including a weighted average in these equations is not required because the decision variable is optimized exclusively. We only need to move the value of the decision variable from the initialization step or previous iteration $k \in K$ at sub-problem $s \in N$ towards the negative direction of the gradient of the objective function. We denote the gradient of the objective function of sub-problem $s \in N$ by $\mathcal{G}_s$ and the step size by $\gamma^{k+1}$. The gradient of objective function of sub-problem $s \in N$ respect to $x_i^{t,od}$ and $y_{ij}^{t,od}$ are shown by $\mathcal{G}_{s,i}^{t,od}$ and $\mathcal{G}_{s,ij}^{t,od}$, respectively.

We use equation (28) for the occupancy variable and (29) for the flow variable to find their auxiliary parameters when a decision variable appears in more than one sub-problem. For instance, the weighted average for decision variable $x_i^{t,od}: t \in T, i \in C^s, (o,d) \in C_{OD}$ at sub-problem $s \in N$ is $\sum_{s,\acute{s} \in n^s, j \in C^{n^s}: i=j} w_{s\acute{s}} x_j^{t,od}$, which takes weighted average of the generated values for this variable by itself and its neighbors from previous iteration $k \in K$. Then, the value of $x_i^{t,od}$ generated by sub-problem $s \in N$ at iteration $k \in K$ is added to term $\alpha^{k+1} \sum_{s,\acute{s} \in n^s, j \in C^{n^s}: i=j} w_{s\acute{s}} x_j^{t,od} - \gamma^{k+1} \mathcal{G}_s$ to find auxiliary parameter $p_i^{t,od}$. We also use $\alpha^{k+1}$ to denote the step size used for consensus among decision variables at iteration $k+1 \in K$. The same approach is used to find other auxiliary parameters at iteration $k + 1 \in K$ as shown in equation (29).

If $x_i^{t,od}: t \in T, i \in C^s, (o,d) \in C_{OD}$ is only in sub-problem $s \in N$: (26)
$$p_i^{t,od} = x_i^{t,od} - \gamma^{k+1} \mathcal{G}_{s,i}^{t,od}$$

If $y_{ij}^{t,od}: t \in T, i \in C^s, j \in S(i), (o,d) \in C_{OD}$ is only in sub-problem $s \in N$: (27)
$$q_{ij}^{t,od} = y_{ij}^{t,od} - \gamma^{k+1} \mathcal{G}_{s,ij}^{t,od}$$

If $x_i^{t,od}: t \in T, i \in C^s, (o,d) \in C_{OD}$ is several sub-problems in addition to $s \in N$: (28)
$$p_i^{t,od} = x_i^{t,od} + \alpha^{k+1} \sum_{s,\acute{s} \in n^s, j \in C^{n^s}: i=j} w_{s\acute{s}} x_j^{t,od} - \gamma^{k+1} \mathcal{G}_{s,i}^{t,od}$$

If $y_{ij}^{t,od}: t \in T, i \in C^s, j \in S(i), (o,d) \in C_{OD}$ is in several sub-problems in addition to $s \in N$: (29)
$$p_i^{t,od} = y_{ij}^{t,od} + \alpha^{k+1} \sum_{s,\acute{s} \in n^s, j \in C^{n^s}: i=j} w_{s\acute{s}} y_{ij}^{t,od} - \gamma^{k+1} \mathcal{G}_{s,ij}^{t,od}$$

Then, we project the value of auxiliary parameters onto the constraints set of sub-problem $s \in N$ as discussed in *Step 3* and using the projection operator described in Definition 4. We find new values for decision variables at iteration $k + 1 \in K$ by solving the following optimization program for sub-problem $s \in N$.

Min $\sum_{(o,d) \in C_{OD}} \sum_{t \in T} \sum_{i \in C^s} \|x_i^{t,od} - p_i^{t,od}\|^2 + \sum_{(o,d) \in C_{OD}} \sum_{t \in T} \sum_{i \in C^s, j \in S(i)} \|y_{ij}^{t,od} - q_{ij}^{t,od}\|^2$ (30)

s.t.
Constraints (2)-(12) for sub-problem $s \in N$

We continue the procedure of updating the value of auxiliary parameters and decision variables for sub-problem $s \in N$ over iterations $k \in K$ until the termination criterion is satisfied.

Equation (31) shows the update procedure in vector notation. The vector of decision variables in sub-problem $s \in N$ at iteration $k + 1$ is denoted by $X_s^{k+1}$, and $X_s$ is the feasible region of sub-problem $s \in N$. Let $w_{s\acute{s}}$ and $\mathcal{G}_s$ denote the weight for link $(s, \acute{s}) \in N$ in information exchange graph $G = (N, A)$ and the gradient of the objective function of sub-problem $s \in N$, respectively. Step sizes at iteration $k + 1 \in K$ are $\alpha^{k+1}$ and $\gamma^{k+1}$.

$$X_s^{k+1} = P_{X_s}[X_s^k + \alpha^{k+1} \sum_{\acute{s} \in n^s} w_{s\acute{s}} X_{\acute{s}}^k - \gamma^{k+1} \mathcal{G}_s]$$ (31)

Fig. 2 shows a small network of four cells that is distributed to two sub-problems. We also illustrate the three steps of the update procedure in Fig. 3 using vector notation to visualize this procedure for this simple example. Each sub-problem contains those constraints and parts of the objective function that have the decision variables corresponding to cells and links within the region assigned to that sub-problem. Note that the constraints and objective function can be distributed following a different structure as long as the explained conditions in section A are satisfied. In this simple example, the sub-problems share the decision variables corresponding to the link between cells 2 and 3. Adding more links between cells 2 and 3 will not change the information exchange graph nor the performance of the approach because this one link forces the flow decision variables to appear in both sub-problems, and these sub-problems share information using two directional arcs in the information exchange graph.

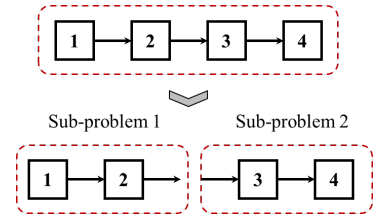

Fig. 2. The distribution of a link with 4 cells to sub-problems 1 and 2

Fig. 3.a shows information exchange graph $G = (N, A)$, where $N = \{1,2\}$ and $A = \{(1,2),(2,1),(1,1),(2,2)\}$. In Fig. 3.b-d, each red circle represents the feasible region of one sub-problem, and vector $z^* \in X^*$ denotes the vector of optimal solutions. Vector $z^* \in X^*$ is within the feasible region of the original problem as well. The number of variables in common among sub-problems will not change this region. Note that sub-problems do not need to share all decision variables. In other words, if there is at least one variable in common, the sub-problems share information.

The vector of initial values for decision variables for sub-problems 1 and 2 at iteration $0 \in K$ are denoted by $X_1^0$ and $X_2^0$,

respectively. Fig. 3.b shows the incorporation of the value of decision variables and the computation of weighted average values for sub-problems 1 and 2 using parameters $a$ and $b$, respectively, as described in *Step 1*. For example, we explain how to compute the value of $a$. We first multiply weight $w_{11}$ on self-arc (1,1) on node 1 in information exchange graph $G$ with the vector of values for decision variables $X_1^0$ in sub-problem 1 at iteration $0 \in K$, that is $w_{11}X_1^0$, and find term $w_{21}X_2^0$ the same way. We then find the value of $a$ by summing terms $w_{11}X_1^0$ and $w_{21}X_2^0$ and multiplying the result by step size $\alpha^1$. The value of $b$ is found following the same procedure. Fig. 3.c presents Step 2 of the update procedure. We first sum the values of $a$ and $b$ with $X_1^0$ and $X_2^0$, respectively. Then, we add $-\gamma^1 G_1$ and $-\gamma^1 G_2$ to the value of $X_1^0 + a$ and $X_2^0 + b$ to move them towards the negative direction of the objective function gradient and find $c$ and $d$, respectively. Fig. 3.d presents the third step, where we project the value of $c$ and $d$ on the feasible region of sub-problems 1 and 2 to find the vector of new values for decision variables $X_1^1$ and $X_2^1$ at iteration $1 \in K$, respectively. In the next section, we explain the required assumption for determining the value of step sizes and weights.

### E. Convergence Properties

It can be shown that the solution of DGA converges to the optimal solution of the SODTA formulation using similar techniques proposed in [47], [48]. Due to space limitations, we only show the assumptions and the final theorem here.

**Assumption 2** The step sizes $\alpha^k$ and $\gamma^k$ should satisfy these conditions: (1) $\sum_{k=1}^{\infty} \alpha^k = \infty$ and $\sum_{k=1}^{\infty} \gamma^k = \infty$, (2) $\sum_{k=1}^{\infty} (\alpha^k)^2 < \infty$ and $\sum_{k=1}^{\infty} (\gamma^k)^2 < \infty$, (3) $\sum_{k=1}^{\infty} (\alpha^k)^2 (\gamma^k)^2 < \infty$, (4) $\sum_{k=1}^{\infty} \frac{(\gamma^k)^2}{\alpha^k} < \infty$, and (5) $\sum_{k=1}^{\infty} \min(\alpha^k, \gamma^k) = \infty$.

**Theorem 1** Vectors $X_s^{k+1}$ generated by sub-problems $s \in N$ using the proposed distributed gradient-based methodology converge to a common optimal vector $z^* \in X^*$, i.e., $\lim_{k \to \infty} X_s^{k+1} = z^*, \forall s \in N$ if Assumption 1 and Assumption 2 hold.

### F. Termination Criterion

The methodology finds the optimal solution when there is no disagreement among the values of decision variables found by different sub-problems. Only some of the flow variables for the boundary links between intersections belong to more than one sub-problem. The optimal solution is found when sub-problems are in agreement on the value of these decision variables, or in other words when the left-hand side of inequality (21) is zero. Since urban streets have a different number of lanes, the value of flow decision variables should be normalized to compute the disagreement by dividing the flow variable by the maximum capacity of receiving or sending cells. We set the termination criterion to reach a disagreement of at most $\varepsilon$.

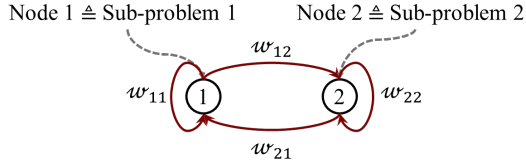

(a) Information exchange graph $G = (N, A)$ such that $N = \{1,2\}$ and $A = \{(1,2),(2,1),(1,1),(2,2)\}$

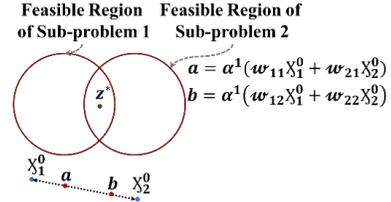

(b) *Step 1*: Finding a weighted average by incorporating the proposed values

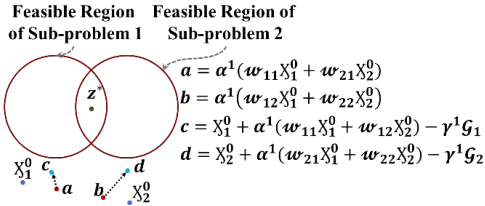

(c) *Step 2*: Moving the values towards the negative direction of the gradient of the objective function from each sub-problem

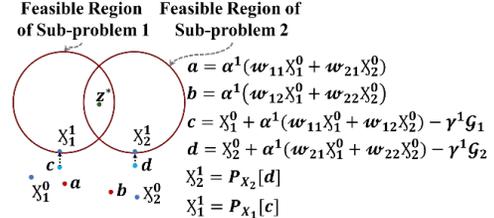

(d) *Step 3*: Projecting the values on the feasible region of each sub-problem

Fig. 3. The procedure of distributed gradient-based update for an iteration

$$\sum_{t \in T} \sum_{od \in OD} \sum_{i \in C^s, j \in S(i), e \in C^{s \in n^s}, f \in S(e): i=e, j=f, s \neq \tilde{s}} \left| \frac{y_{ij}^{t,od}}{\max(M_i, M_j)} - \frac{y_{ef}^{t,od}}{\max(M_e, M_f)} \right| \leq \varepsilon \quad (32)$$

This termination criterion is checked using a distributed communication paradigm. Each sub-problem checks inequality (32), and if it satisfies the termination criterion, it stops updating the value of variables. The sub-problems compute the termination at the same time.

## V. TEST NETWORK

We tested DGA on a Nguyen and Dupuis [49] network and portion of the downtown Springfield network in Illinois. We used the modified Nguyen and Dupuis [49] network by Long and Szeto [21]. This network contains 17 nodes, 23 links, and four OD pairs, and its cell representation has 57 cells, 63 links, 4 OD pairs. The network is decomposed into two regions creating two sub-problems in DGA. The two regions have 30 and 27 cells. They also have 36 and 31 links, respectively, from which 4 links are shared and the average number of decision variables in each subproblem is 2700 variables. Similar demand level, jam density, and flow capacity in [20] have been used. The network consists of 20 intersections with one-way and two-way streets, as shown in Fig. 4. All intersections are signalized with predefined signal timing parameters. In this figure, roads are shown with links, and the circles show the intersection number.

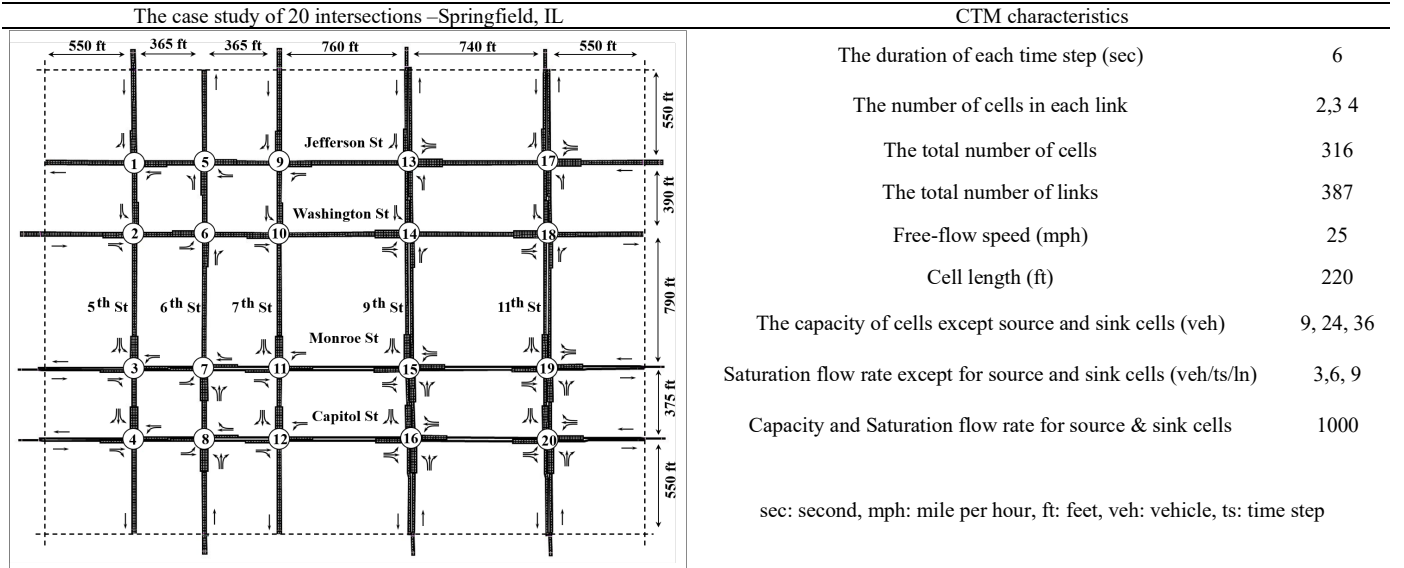

Fig. 4. The case study of 20 intersections, its CTM characteristics, and demand patterns

We considered 15 OD pairs and three demand profiles for this test network, representing undersaturated, saturated, and oversaturated conditions. We assumed demand for ODs is known. The proposed technique in [50], [51] may be used to estimate and calibrate the demand in the future. A test network with 40 (4×10) intersections and 25 ODs is also tested to show how the approach scales. This network is created by duplicating the network of 20 intersections. We used three demand profiles to represent different saturation levels. According to Assumption 1, one set of values for weights can be $w_{s\acute{s}} = w_{\acute{s}s} = 1$ for all $s, \acute{s} \in N$ in which $s \neq \acute{s}$ and $w_{ss} = -\sum_{\acute{s} \in C^{ns}: s \neq \acute{s}} w_{s\acute{s}}$ for weights on self-arcs $s \in N$. We also set the step size rules as $\alpha^k = 1/k^{0.55}$ and $\gamma^k = 1/k$, that satisfy Assumption 2. For 20 intersection network, the termination is set to reach a disagreement of 0.5 for each sub-problem. It should be noted that networks with 20 and 40 intersections respectively result in 4,218,000 and 14,120,000 decision variables, representing a large optimization program.

VI. RESULTS

A. Benchmarking the Methodology

- Comparison with a Central Approach

We applied DGA and a central approach to the case study network of 20 intersections. TABLE II shows the total travel time and the total computation time by our approach and the central approach under three different demand patterns. The SODTA formulation has more than 4 million decision variables in this case study, and we could not find the optimal solution by running CPLEX for cases with less than 150 GB of memory. DGA reduced the decision variables of each sub-problem to 243,600 on average. Each pair of sub-problems shares 6,000 variables for one-directional street and 12,000 variables for two-directional streets. DGA required only 5 GB of memory to generate the solutions. Note that the distributed approach is implemented with a parallel architecture. It also found the solutions with at most 5% optimality gap in less than 2.01 hours, which translates to about 97 % shorter runtime than CPLEX.

Moreover, the computation time of the algorithm proposed in this paper is 0.76 seconds, which is 33.34% less than the lowest CPU time for all tested formulations by Long et al. [20]. The number of variables for LTM and CTM-based formulations are 1,656 and 5,400. The CPU time to solve four different CTM-based formulations by Long and Szeto [21] varies between 1.14 seconds to 36.32 seconds.

TABLE II Objective function and computation time for 20-intersection NetWOrk

| Approach/Gap | Demand Profile | Memory (GB) | Objective Function: Total Travel Time (hr) | Run-time (hr) |
|---|---|---|---|---|
| CPLEX | 1 | 150 | 187.06 | 70.54 |
| DGA | 1 | 5 | 197.34 | 1.60 |
| Difference (%) | 1 | - | 5.49 | 97.73 |
| CPLEX | 2 | 150 | 291.29 | 73.59 |
| DGA | 2 | 5 | 307.30 | 1.81 |
| Difference (%) | 2 | - | 5.50 | 97.54 |
| CPLEX | 3 | 150 | 499.85 | 90.56 |
| DGA | 3 | 5 | 524.84 | 2.01 |
| Difference (%) | 3 | - | 5.01 | 97.78 |

Fig. 5.a-c shows the objective value of DGA and the optimal objective value for three demand profiles over iterations. The distributed approach reduces the value of the objective function towards the optimal value over iterations. Note that the optimal solution from CPLEX had the well-known holding back problem due to the linearization of the minimization functions in CTM. However, the solutions from the distributed approach did not have this issue since we found the solutions by simulating the network using non-linear CTM equations as presented in [52]. In all three figures, the optimality gap was at most 5% when we stopped the approach.

- *Comparison with a Danzig-Wolfe Decomposition-based Algorithm (DWDA)*

We compared the solutions and performance of the DGA to an OD-based decomposition approach developed in [14]. TABLE III shows different characteristics and performance measures for both approaches. Danzig-Wolfe Decomposition-based Algorithm (DWDA) has a master problem and several sub-problems that are solved iteratively with a stopping criterion of a 5% gap between the upper bound and lower bound of the approach. Each sub-problem has all decision variables and constraints for a SODTA formulation with one OD pair. Since we have 15 ODs, the number of sub-problems is 15. New extreme points are generated by sub-problems and added to the master problem solution pool. Therefore, the complexity of the master problem increases over iterations.

total run-time of DGA is 74% more than DWDA in undereducated demand pattern. However, when most of the decision variables have non-zero values in the oversaturated condition, the runtime of DGA is improved by 77% compared to DWDA. Moreover, since we have 48% fewer variables in DGA, we only required 5 GB of memory though we need at least 20 GB of memory to run DWDA.

TABLE III Benchmark with Danzig-Wolfe Decomposition algorithm for the network of 20 intersections

| Demand | Approach/ Gap | Number of variables in sub-problems | Iterations | Objective Function | Optimality gap (%) | Run-time (hr) | Memory (GB) |
|---|---|---|---|---|---|---|---|
| 1 | DWDA | 281,200 | 11 | 197.69 | 0.05 | 0.92 | 20 |
| 1 | DGA | ≤ 146,880 | 1706 | 197.34 | 0.05 | 1.60 | 5 |
| 1 | Diff. (%) | -48 | 15409 | 0.18 | 0.00 | 74 | -75 |
| 2 | DWDA | 281,200 | 53 | 308.42 | 0.05 | 2.34 | 20 |
| 2 | DGA | ≤ 146,880 | 1708 | 307.30 | 0.05 | 1.81 | 5 |
| 2 | Diff. (%) | -48 | 3123 | 0.36 | 0.00 | -23 | -75 |
| 3 | DWDA | 281,200 | 101 | 529.29 | 0.05 | 8.83 | 20 |
| 3 | DGA | ≤ 146,880 | 1715 | 524.84 | 0.05 | 2.01 | 5 |
| 3 | Diff. (%) | -48 | 1599 | 0.84 | 0.00 | -77 | -75 |

### B. The Performance of the Methodology

Fig. 6 presents the disagreement on the value of shared decision variables for three sub-problems with their neighbors for three demand profiles over iterations. In demand profile 1, when each sub-problem has a disagreement value of less than 0.5, the algorithm is terminated.

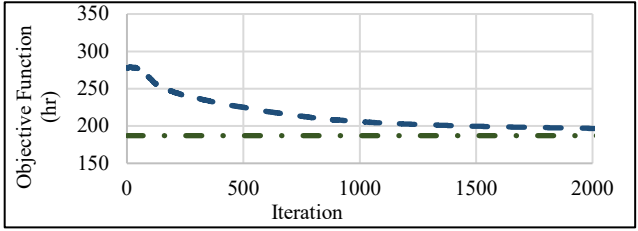
(a) demand profile 1

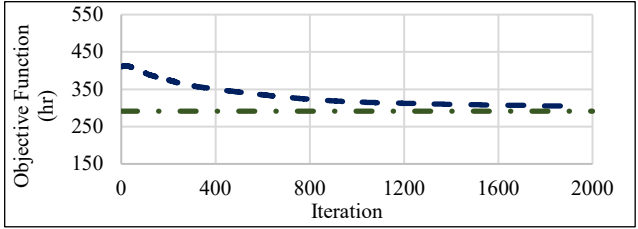
(b) demand profile 2

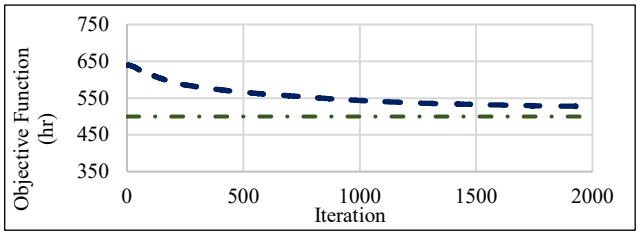
(c) demand profile 3

– – Distributed Algorithm  — · · Optimal Solution

Fig. 5. The objective functions for the network of 20 intersections

DGA has an intersection-based decomposition, and the number of sub-problems is 20 due to having 20 intersections in the network. This approach does not have any central component or master problem. The number of decision variables differs slightly in each sub-problem of DGA depending on the number of nodes and links and is at least 48% less than the number of variables in the sub-problems of DWDA. The number of iterations in DWDA is at most 101 while this number is 1715 for DGA. Even though we have more number of iterations in DGA, the computation time of each iteration is much less. The

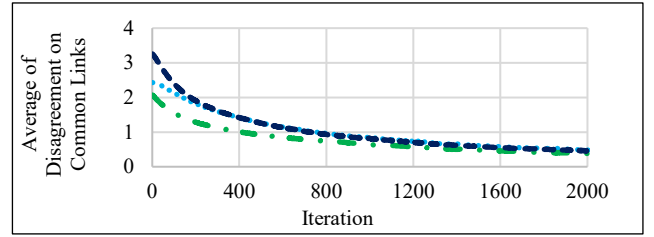
(a) demand profile 1

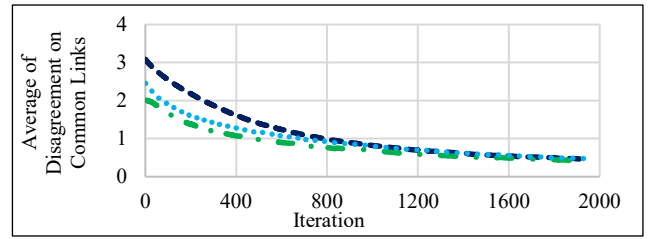
(b) demand profile 2

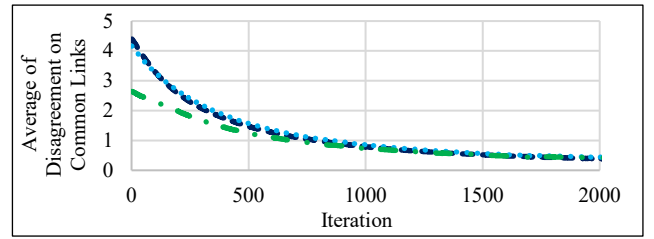
(c) demand profile 3

– – – Subproblem 1   ······ Subproblem 2   — · — Subproblem 3

Fig. 6. Disagreement value $\varepsilon$ over iterations in 20-intersection network

The algorithm reached a 5.5% optimality gap for demand profile 1 at iteration 1706, see Fig. 6.a. The algorithm reached a disagreement value of 0.5 for the second demand pattern at iteration 1708 with the 5.5% optimality gap, see Fig. 6.b. The approach also reached a disagreement value of 0.5 at iteration 1715 iterations for the third demand pattern with a similar gap of 5.0%. As expected, the disagreement increased with the demand level. However, the gap is the same for all scenarios once the approach is terminated.

Fig. 7.a - c shows the run-time of three sub-problems over iterations for the three demand profiles. We assigned each sub-problem to a different computational node using a multi-thread platform. The sub-problems are independent and optimized synchronously. Each sub-problem represents one intersection with an almost equal number of variables and constraints with other sub-problems. Therefore, the run-time is approximately the same among different nodes, which reduces overhead delays. The run-time of each sub-problem is relatively similar for all three demand profiles, and having a more congested sub-problem did not significantly affect the computation time; however, the total run-time differs by at most 20% between demand profiles due to the additional number of required iterations.

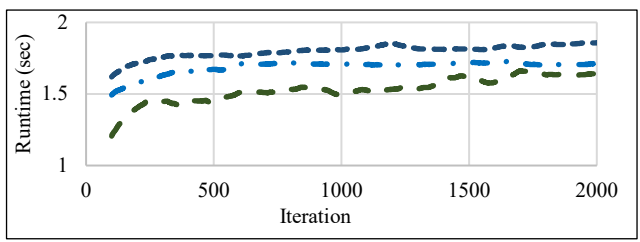

(a) demand profile 1

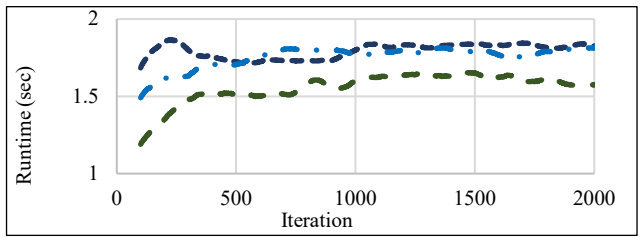

(b) demand profile 2

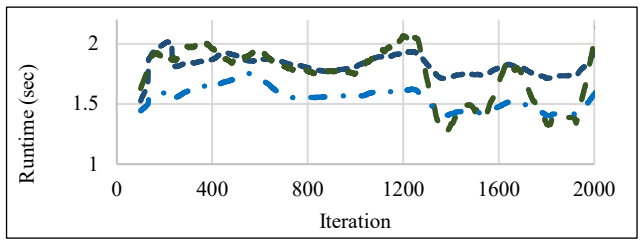

(c) demand profile 3

----- Subproblem 1   — — Subproblem 2   - · - Subproblem 3

Fig. 7. The run-time of three sub-problems 20-intersection network

TABLE IV shows the computation time for each section of the approach under three demand patterns. The runtime for the initialization step is negligible for all demand profiles. The distributed gradient update consists of step 1: computing weighted averages, step 2: moving values towards the negative direction of the gradient, and step 3: projecting values on the feasible region of sub-problems. Finding the weighted average needs the least CPU time, and projection has the highest computation time. The projection of the value of variables is implemented for all sub-problems at the same time, and its CPU time has increased by increasing the demand. This can happen because more number of decision variables has non zero values as the demand increases, and more computational effort is required to find variables' optimal values. No specific trend has appeared for other steps including termination criterion calculation.

TABLE IV
Breakdown of runtimes for the network of 20 intersections

| Demand Level | | 1 | 2 | 3 |
|---|---|---|---|---|
| Initialization (hr) | | 0.0000 | 0.0000 | 0.0000 |
| Distributed gradient update (hr) | Step 1 | 0.0023 | 0.0024 | 0.0022 |
| | Step 2 | 0.0035 | 0.0036 | 0.0033 |
| | Step 3 | 1.5525 | 1.7643 | 1.9689 |
| Termination criterion (hr) | | 0.0351 | 0.0346 | 0.0338 |

Fig. 8.a-d shows the disagreement on the value of decision variables for each sub-problem with its neighbors at iterations 1, 500, 1000, 1800 for three demand profiles. Each sub-problem is shown with a number following the same layout shown in Fig. 4. The spectrum shows a range of colors depending on the value of the disagreement. Darker color represents a higher disagreement value. Fig. 8.a-d present the disagreement at each sub-problem for the under-saturated demand profile. The colors become lighter as the number of iterations increases, which shows that the conflict on the proposed value by each sub-problem with its neighbors decreases. Fig. 8.e-h and Fig. 8.j-m display the same pattern for the second and third demand profiles, respectively.

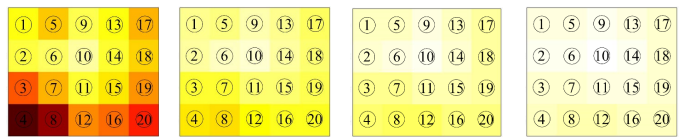

(a) Iteration 1    (b) Iteration 500    (c) Iteration 1000    (d) Iteration 1800
Demand profile 1

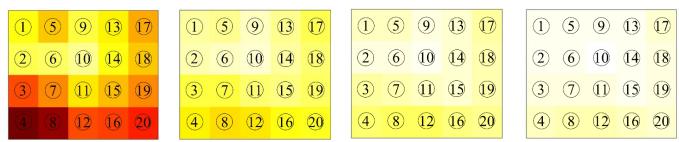

(e) Iteration 1    (f) Iteration 500    (g) Iteration 1000    (h) Iteration 1800
Demand profile 2

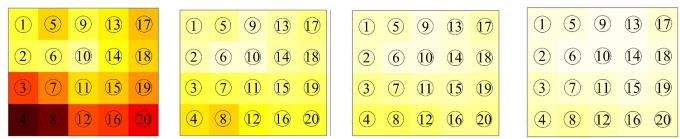

(e) Iteration 1    (f) Iteration 500    (g) Iteration 1000    (h) Iteration 1800
Demand profile 3

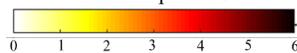

Fig. 8. The value of disagreement at each sub-problem.

## C. The Effects of the Number of ODs on Convergence

Fig. 9 shows the impact of the number of OD pairs (15, 20, and 40) on the convergence of our approach in the test network of 20 intersections and semi-saturated demand. The algorithm reached the termination criteria in 1708, 2157, and 2319 iterations with 5%, 4%, and 5% optimality gaps for 15, 20, and 40 ODs, respectively. Increasing the number of OD pairs led to more iterations for convergence; however, the number of iterations does not increase as fast as the number of decision variables. Specifically, increasing the number of ODs from 20 to 40 doubles the number of decision variables but only increases the number of iterations by 7%.

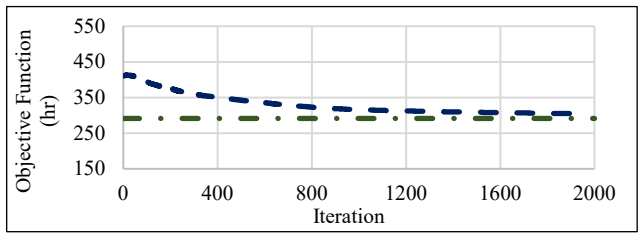

(a) 15 OD pairs

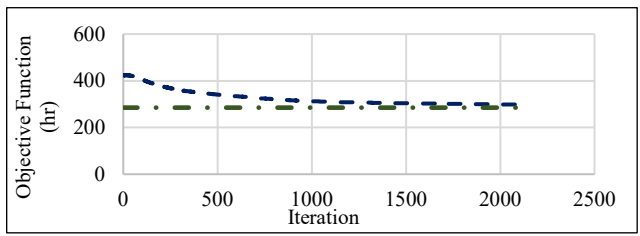

(b) 20 OD pairs

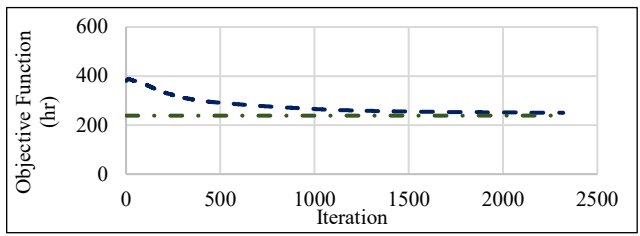

(c) 40 OD pairs

Fig. 9. The effects of the number of OD pairs on DGA convergence

## D. The Effects of the Study Period Duration on Convergence

In this section, we evaluated the effect of increasing the duration of the study period on convergence. We tested three loading periods of 150, 300, and 450 time steps, as shown in Fig. 10. The three cases required 1420, 1708, and 2081 iterations with 5%, 4%, and 4% optimality gaps to reach the termination criterion. As we increased the loading time, the number of required iterations to meet the termination was increased. Increasing the loading time created more congestion. Therefore, the number of iterations was increased by increasing the loading period. However, the rate of increase in the required iterations was much less than the increase in the number of decision variables.

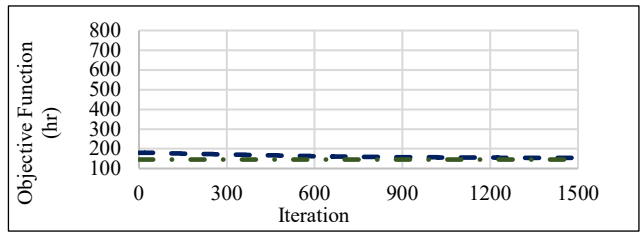

(a) 150 time steps for loading

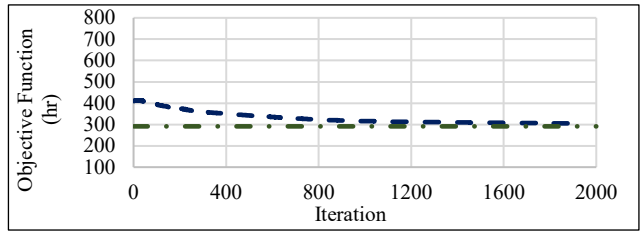

(b) 300 time steps for loading

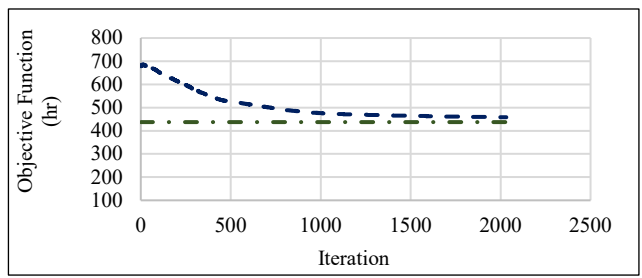

(c) 450 time steps for loading

Fig. 10. The effect of study period duration on DGA convergence

## E. The Effects of the Network Size on Convergence

We also studied the effects of network size on convergence by looking at networks of 10, 20, and 40 intersections with similar characteristics, as shown in Fig. 11. An increase in network size increased the number of decision variables from 1,674,000 to 4,218,000, and 14,120,00. DGA was converged in 807, 1708, and 1710 iterations with 2%, 4%, and 5% optimality gaps for networks with 10, 20, 40 intersections, respectively. The loading period, demand, and the number of OD pairs were the same in all cases. The number of iterations for convergence was increased by increasing the network size; however, at a rate much slower than the increase in the network size. Increasing the size of the network from 20 intersections to 40 increases the number of required iterations by only two, which shows the scalability of the proposed methodology.

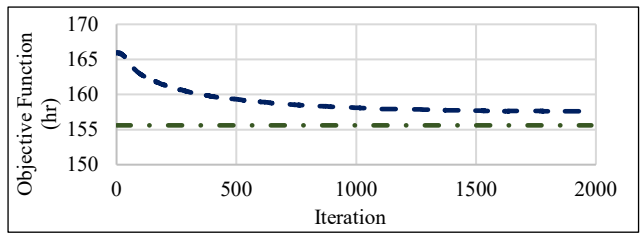

(a) 10 intersections

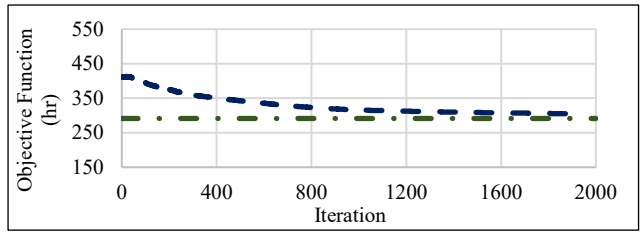

(b) 20 intersections

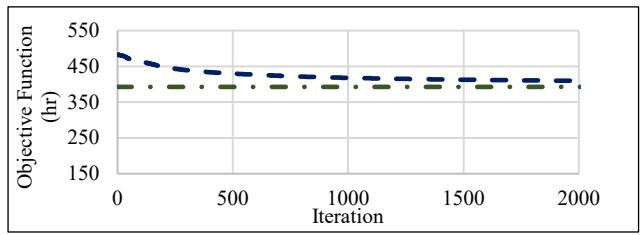

(c) 40 intersections

— — Distributed Algorithm  — · — Optimal Solution

Fig. 11. The effect of network size on DGA convergence

*F. The effect of Network Size on Decision Space*

We also tested the approach on a network with 40 intersections with 632 cells, 780 links, 400 time steps, and 25 OD pairs that bring the total number of decision variables to 14,120,000. We used the same termination criterion that was applied to other cases. By increasing the intersections from 20 to 40, the number of decision variables increased from ~4 million to ~14 million (more than a factor of 3). We need to mention that the literature that includes models of networks with thousands of intersections and OD pairs use either exit flow functions, point queue models, or link-performance functions. These approaches are aggregated and have significantly fewer decision variables and do not provide the accuracy that is required for traffic operation purposes. In this paper, we use the CTM model, which is more accurate but at the expense of additional complexity. The network of 40 intersections is significantly larger than comparable studies that have used the cell transmission model [11], [13], [53]. The studies solve the problem for 5805 to 489,700 decision variables [12], [17].

Increasing the number of links from 387 to 780 and nodes from 316 to 632 does not change the complexity of each sub-problem due to the intersection-based distribution of the formulation. However, increasing the study period from 400 to 500 time steps and OD pairs from 15 to 25 lead to having more variables. Even though the variables are increased by 70%, we can find the solutions with at most 5.70% optimality gap in at most 5.34 hours, as shown in TABLE V.

TABLE V increasing the network size from 20 to 40 intersections

| Demand | Intersections | Optimality gap (%) | | Runtime (hr) | | Number of Iterations | |
|---|---|---|---|---|---|---|---|
| | | DGA | DWDA | DGA | DWDA | DGA | DWDA |
| 1 | 20 | 5.90 | 5.06 | 1.60 | 0.92 | 1704 | 11 |
| | 40 | 5.68 | NA | 5.19 | NA | 1710 | NA |
| 2 | 20 | 5.28 | 5.88 | 1.81 | 2.34 | 1708 | 53 |
| | 40 | 5.68 | NA | 5.26 | NA | 1711 | NA |
| 3 | 20 | 5.93 | 5.91 | 2.01 | 8.83 | 1715 | 101 |
| | 40 | 5.70 | NA | 5.34 | NA | 1716 | NA |

*G. The effect of the Network Size on Run-time*

Fig. 12 shows the total run-time for sub-problems for the network of 20 intersections with 15 and 30 ODs and the network of 40 intersections with 25 ODs. Increasing the number of OD pairs from 15 to 30 in the network of 20 intersections has increased the number of decision variables by 50% and run-time at each iteration by 71% on average. When the size increases, the number of variables in each sub-problem increases only by the number of OD pairs in this case. The number of cells, links, and time steps remains constant. By increasing the number of intersections from 20 with 15 ODs to 40 intersections with 25 ODs, the number of decision variables in each sub-problem has increased by 70%. The run-time to generate solutions at each time step has increased by 52% on average.

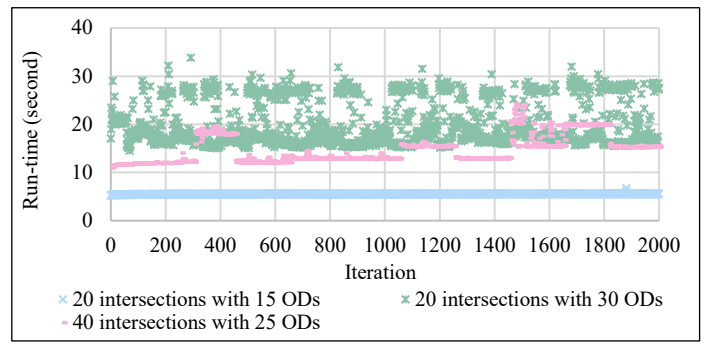

Fig. 12. The run-times of the proposed methodology

## VII. CONCLUSION

We introduced a gradient-based methodology that distributed SODTA into several independent sub-problems. This methodology overcomes the main drawback of existing decomposition techniques by having a fully distributed framework that does not require a centralized component. Each sub-problem represents one intersection. The set of constraints and the objective function of the CTM-based SODTA problem is distributed among sub-problems such that each sub-problem has some parts of the SODTA formulation with decision variables corresponding to the intersection region assigned to that sub-problem.

DGA updates the decision variables at each iteration by combining the proposed values for the shared decision variables by finding their weighted average. The approach also considers the gradient of the objective function at each sub-problem in updating these values. Then, DGA projects the computed values on the feasible region of each sub-problem. The update procedure continues until it reaches a consensus on the proposed values for shared decision variables among sub-problems. DGA was tested on a Nguyen and Dupuis [49] network and a portion of the downtown Springfield network in Illinois. The run time of this approach was 33.34% less than the lowest run time for tested formulations on Nguyen and Dupuis [49] network by Long and Szeto [21] . DGA required 96% less memory and 97% less time compared to CPLEX to find solutions with a 5% optimality gap. In comparing DGA with a Danzig-Wolfe Decomposition-based Algorithm, the runtime was improved by 77% in oversaturated condition.

Adapting DGA to solve various types of DTA formulations like the DTA model proposed for emergency conditions [54] and the DTA model incorporated with signal timing optimization [55]–[57] can be considered as future research directions. Further research on different network partitioning scenarios presented in [58] and their effects on the convergence and run-time of our approach is needed. In a future study, we will consider multi-modal transportation systems similar to the model proposed by Bevrani et al. [59]. A more general formulation with multi-commodity flows similar to the presented model by Bevrani et al. [60], [61] can also be employed for generalizing this approach. We will also study the possibility of estimating the expected number of iterations to achieve a solution within a small gap in future studies.

APPENDIX
TABLE VI Notations

| | |
|---|---|
| **Sets** | |
| $T$ | The set of time steps |
| $C$ | The set of network cells |
| $C_O$ | The set of source cells |
| $C_S$ | The set of sink cells |
| $C_I$ | The set of intersection cells |
| $C_{OD}$ | The set of OD pairs |
| $P(i)$ | The set of predecessors of cell $i \in C$ |
| $S(i)$ | The set of successors of cell $i \in C$ |
| $C_D$ | The set of diverge cells |
| $N$ | The set of sub-problems |
| $Cl$ | The set of links that their heads $i \in C^s$ and tails $j \in S(i)$ belong to sub-problems $e \in C^{s \in n^s}$ and $f \in S(e)$ such that $i = e, j = f, s \neq s$ |
| $\mathfrak{n}^s$ | The set of neighbors of sub-problem $s \in N$ including itself |
| $C^s$ | The set of cells in sub-problem $s \in N$ |
| $Co^s$ | The union of the set of cells in sub-problem $s \in N$ and the set of source cells, i.e., $C_O \cap C^s$ |
| $Cs^s$ | The union of the set of cells in sub-problem $s \in N$ and the set of sink cells, i.e., $C_S \cap C^s$ |
| $C_D^s$ | The set of diverge cells in sub-problem $s \in N$ |
| $Ceos^s$ | The set of cells in sub-problem $s \in N$ except for the set of source cells and the set of sink cells i.e., $C^s/C_O \cup C_S$ |
| $Ces^s$ | The set of cells in sub-problem $s \in N$ except for the set of sink cells, i.e., $C^s/C_S$ |
| $Ceo^s$ | The set of cells in sub-problem $s \in N$ except for the set of source cells, i.e., $C^s/C_O$ |
| $Ci^s$ | The union of the set of cells in sub-problem $s \in N$ and the set of intersection cells, i.e., $C_I \cap C^s$ |
| $X_s$ | The set of constraints of subproblem $s \in N$ |
| $X$ | The feasible region of SODTA problem |
| $X^*$ | The set of optimal solutions |
| **Decision variables** | |
| $x_i^{t,od}$ | The number of vehicles in cell $i \in C$ at time step $t \in T$ for OD pair $(o,d) \in C_{OD}$ |
| $y_{ij}^{t,od}$ | The number of vehicles flowing from cell $i \in C$ to downstream cell $j \in S(i)$ at time step $t \in T$ for OD pair $(o,d) \in C_{OD}$ |
| **Parameters** | |
| $\tau$ | The duration of time step $t \in T$ |
| $\mathcal{D}_i^{t,od}$ | The entry demand level at source cell $i \in C_O$ at time step $t \in T$ for OD pair $(o,d) \in C_{OD}$ |
| $F_i$ | The saturation flow rate at cell $i \in C$ |
| $M_i$ | The maximum number of vehicles that cell $i \in C$ can accommodate |
| $\delta$ | The ratio of free-flow speed to the backward propagation speed |
| $g_i^t$ | A binary parameter to define signal status at intersection cell $i \in C_I$ at time step $t \in T$. Zero and one values indicate red and green signals, respectively |
| $R_{ij}^{t,od}$ | The turning ratio of the link between diverge cell $i \in C_D$ and its successor cell $j \in S(i)$ at time step $t \in T$ for OD pair $(o,d) \in C_{OD}$ |
| $f_i^t$ | The variable saturation flow rate of intersection cell $i \in C_I$ at time step $t \in T$ |
| **Notations and Terms** | |
| $K$ | The total number of iterations |
| $k$ | The iteration counter |
| $X_s^k$ | The column vector of values for the decision variables of sub-problem $s \in N$ at iteration $k \in K$ |
| $x$ and $z$ | The column vectors |
| $x'$ | The transpose of vector $x$ |
| $x'z$ | The dot product of two vectors $x$ and $z$ |
| $z^* \in X^*$ | A realization of the vector of optimal solution |
| $\|x\| = (x'x)^{\frac{1}{2}}$ | The standard Euclidean norm |
| $P_X[x]$ | The projection of vector $x$ on set $X$ |
| $\mathcal{F}$ | The objective function of SODTA problem |
| $\mathcal{F}_s(x)$ | The objective function of sub-problem $s \in N$ at vector $x$ |
| $\mathcal{F}^*$ | The optimal objective function |
| $G = (N, A)$ | The information exchange graph with a set of nodes $N$ and a set of directional links $A$ |
| $(s, \mathfrak{s}) \in A$ | A link in information exchange graph $G = (N, A)$ |
| $w_{s\mathfrak{s}}$ | The weight on the link going from node $s$ to node $\mathfrak{s}$ in information exchange graph $G = (N, A)$ |
| $\theta$ and $\theta'$ | The lower bound and the upper bound for weight values $w_{s\mathfrak{s}}$: $s\mathfrak{s} \in N$ |
| $w$ | The weighted graph Laplacian matrix with weight $w_{s\mathfrak{s}}$ entries assigned to links $(s, \mathfrak{s}) \in N$ in information exchange graph $G = (N, A)$ |
| $\mathcal{G}_s$ | The gradient of the objective function $\mathcal{F}_s$ of sub-problem $s \in N$ |
| $\mathcal{G}_{s,i}^{t,od}$ | The gradient of objective function of sub-problem $s \in N$ respect to $x_i^{t,od}$ |
| $\mathcal{G}_{s,ij}^{t,od}$ | The gradient of objective function of sub-problem $s \in N$ respect to $y_{ij}^{t,od}$ |
| $\alpha^k$ and $\gamma^k$ | The step sizes at iteration $k \in K$ |
| $\varepsilon$ | The consensus error used for termination criterion |
| $p_i^{t,od}$ | The auxiliary parameter for vehicles in cell $i \in C$ at time step $t \in T$ with OD pair $(o,d) \in C_{OD}$ |
| $q_{ij}^{t,od}$ | The auxiliary parameter for vehicles going from cell $i \in C$ to downstream cell $j \in S(i)$ at time step $t \in T$ with OD pair $(o,d) \in C_{OD}$ |
| $\mathbb{A} \setminus \mathbb{B}$ | All elements in set $\mathbb{A}$ except for the ones in set $\mathbb{B}$ |